\RequirePackage[l2tabu, orthodox]{nag}
\documentclass[preprint,  11pt]{amsart}

\usepackage{fullpage}

\usepackage{amsfonts}
\usepackage{amssymb}
\usepackage{amsmath}
\usepackage{amsthm}
\usepackage{amsbsy}
\usepackage{amssymb} 
\usepackage{verbatim}
\usepackage{bm} 
\usepackage{paralist} 
 \usepackage{color}
 \usepackage{mathrsfs}
 \usepackage{graphicx}
\usepackage{hyperref}

\newcommand{\E}{\mathbf{E}}
\def\P{\mathbf{P}}

\def\given{\left.\vphantom{\hbox{\Large (}}\right|}
\def\Given{\ \pmb{\big|} \ }

\newcommand{\mat}[4]{\l[\begin{array}{cc}  #1 & #2 \\ #3 & #4  \end{array} \r]}

\newcommand{\one}{{\mathbf 1}}

\def\summ{\sum\limits}
\def\intt{\int\limits}

\def\prodd{\prod\limits}

\def\mb{\mbox}

\newcommand{\para}[1]{\vspace{4mm}\noindent{\bfseries #1:}}

\def\l{\left}
\def\r{\right}
\def\<{\langle}
\def\>{\rangle}

\newcommand{\ba}{\[\begin{aligned}}
\newcommand{\ea}{\end{aligned}\]}
\newcommand\mnote[1]{} 
\newcommand{\beq}[1]{\begin{equation}\label{#1}}
\newcommand\eeq{\end{equation}}
\newcommand\ben{\begin{equation}}
\newcommand\een{\end{equation}}
\newcommand\bes{\begin{eqnarray*}}
\newcommand\ees{\end{eqnarray*}}
\newcommand\besn{\begin{eqnarray}}
\newcommand\eesn{\end{eqnarray}}

\def\bthm{\begin{theorem}}
\def\ethm{\end{theorem}}
\def\bdefn{\begin{definition}}
\def\edefn{\end{definition}}
\newcommand{\benu}{\begin{enumerate}\setlength\itemsep{6pt}}
\newcommand{\beit}{\begin{itemize}\setlength\itemsep{3pt}}
\def\eenu{\end{enumerate}}
\def\eeit{\end{itemize}}
\def\beds{\begin{description}}
\def\eeds{\end{description}}
\def\bepr{\begin{problem}}
\def\eepr{\end{problem}}

\def\bprf{\begin{proof}}
\def\eprf{\end{proof}}
\def\berk{\begin{remark}}
\def\eerk{\end{remark}}
\def\bex{\begin{exercise}}
\def\eex{\end{exercise}}
\def\beg{\begin{example}}
\def\eeg{\end{example}}

\def\suchthat{{\; : \;}}
\renewcommand{\qed}{\hfill\text{$\blacksquare$}}

\def\R{\mathbb{R}}

\def\Z{\mathbb{Z}}

\newcommand{\sm}{{\raise0.3ex\hbox{$\scriptstyle \setminus$}}}

\def\alp{\alpha}
\def\bet{\beta}
\def\gam{\gamma}
\def\del{\delta}
\def\eps{\epsilon}
\renewcommand\phi{\varphi}

\def\lam{\lambda}

\newcommand\x{\mathbf{x}}

\theoremstyle{plain} 
    \newtheorem{theorem}{Theorem}
    \newtheorem{lemma}[theorem]{Lemma}
    
    \newtheorem{corollary}[theorem]{Corollary}

\theoremstyle{definition} 
    \newtheorem{definition}[theorem]{Definition}

       \newtheorem{problem}[theorem]{\bf }

    \newtheorem{exercise}[theorem]{Exercise}
        \newtheorem{remark}[theorem]{Remark}
    \newtheorem{example}[theorem]{Example}

\renewcommand{\bex}{\indent\begin{exercise}}

\usepackage{palatino}
\usepackage{undertilde}

\hypersetup{
    colorlinks=true, 
    linktoc=all,     
    linkcolor=blue,  
}
\begin{document}

\title{A relative anti-concentration inequality}
\author{Manjunath Krishnapur and Sourav Sarkar}
\address{Department of Mathematics\\
        Indian Institute of Science\\
        Bangalore 560012, India}
\email{manju@math.iisc.ernet.in}
\thanks{}

\address{Department of Statistics, University of California, Berkeley, USA}
\email{souravs@berkeley.edu}
\thanks{The second author was supported in part by IAS summer research fellowship and Lo\`{e}ve Fellowship at University of California, Berkeley.}

\begin{abstract} Given two vectors in Euclidean space, how unlikely is it that a random vector has a larger inner product with the shorter  vector than with the longer one? When the random vector has independent, identically distributed components, we conjecture that this probability is no more than a constant multiple of the ratio of the Euclidean norms of the two given vectors, up to an additive term  to allow for the possibility that the longer vector has more arithmetic structure.  We give some partial results to support the basic conjecture. 
\end{abstract}

\maketitle

\section{The question}
We conjecture the following {\em relative anti-concentration} inequality: If  $\alp,\bet\in \R^{n}$, and $X_{i}$ are  i.i.d. real-valued random variables with a non-degenerate distribution, then
\begin{align}\label{eq:relanticoncentrationgoal}
\P\l\{\Big| \sum_{i=1}^{n}\alp_{i}X_{i} \Big|\le \Big|\sum_{i=1}^{n}\bet_{i}X_{i} \Big|\r\}\le C\frac{\|\bet\|}{\|\alp\|}+\frac{C}{\mb{LCD}(\alp)}.
\end{align}
Here $C$ is a constant, $\|\beta\|^{2}=\bet_{1}^{2}+\ldots +\bet_{n}^{2}$, and $\mb{LCD}(\alp)$ is the  ``essential least common denominator'' introduced by Rudelson and Vershynin in their inverse Littlewood-Offord theorems. Its precise definition is recalled later.  In this paper, we prove special cases of this inequality, under conditions on the distribution of $X_{1}$ or on the coefficients, and in some cases not requiring the second term at all. 

To put the inequality in context, recall the L\'{e}vy concentration function of a real-valued  random variable $X$, defined as $$Q_{X}(t)=\sup_{a\in \R}\P\{a\le X\le a+t\}.$$ {\em Anti-concentration inequalities} are upper bounds on the concentration function, perhaps for a  range of $t$ (for instance, on $Q_{X}(0)$, which is the maximal size of an atom). The famous Littlewood-Offord problem is an anti-concentration inequality for $S=\sum_{i=1}^{n}v_{i}X_{i}$ where $X_{i}$ are independent Bernoulli random variables. It states that $Q_{S}(t)\le C/\sqrt{n}$, provided $t\le v_{i}$ for all $i$. This has been  generalized in different directions. The Kolmogorov-Rogozin inequality generalizes to  sums of independent random variables. H\'{a}lasz's inequalities and the inverse Littlewood-Offord theorems (Arak, Tao and Vu, Rudelson and Vershynin, etc.) are stronger bounds on $Q_{S}$ (also allowing general distributions of $X_{i}$s) under constraints on the arithmetic structure of $v_{i}$s.  See  \cite{nguyenvu}, \cite{gotzezaitsev} or \cite{rudelsonvershyninicm} for more on this fascinating subject. In short, these are upper bounds on the small-ball probabilities of linear forms under product measure. There are  anti-concentration theorems of quadratic forms of independent random variables and more generally for polynomials (eg., \cite{taovubook}, \cite{mekanguyenvu}).

Now it is clear why we call \eqref{eq:relanticoncentrationgoal} a ``relative'' anti-concentration inequality (think of $\|\beta\|$ as small and $\|\alp\|$ as large, else the inequality is trivial), since it asks for the probability that a linear form with small coefficients dominates another one with large coefficients. Why do we expect the bound on the right? If  $X_{i}$ are i.i.d. standard Gaussians, then it is an easy calculation (shown later) that the probability is bounded by $C\|\bet\|/\|\alp\|$. We expect essentially the same bound in general, except that for discrete random variables such as Bernoullis, the second term is needed. This is because the quantity $\|\beta\|/\|\alp\|$ can be made as small as desired by scaling $\beta$ down, while the left hand side cannot be smaller than the atom size of $\sum \alp_{i}X_{i}$ at $0$ (which can be non-zero if $\alp$ has an arithmetic structure). The term $1/\mb{LCD}(\alp)$ is precisely what Rudelson and Vershynin use to bound the largest atom of $\sum \alp_{i}X_{i}$.

The special case when $\beta_{i}=1$  and $\alp_{i}=i$, has an application to the study of zeros of random polynomials. In this case, the inequality \eqref{eq:relanticoncentrationgoal} (the bound on the right is simply $1/n$) was proved by S\"{o}ze~\cite{kensoze} (see Lemma~3 in his paper) who  used it to prove a bound for the expected number of real zeros of random polynomials with i.i.d. coefficients. Other than that, we do not know of any applications of the inequality~\eqref{eq:relanticoncentrationgoal}. However it appears to have a natural appeal and in this paper we prove several partial results to support our conjecture.

\para{Acknowledgement} After the first version of our article was posted on the axiv, Sasha Sodin communicated to us a Fourier analytic proof of \eqref{eq:relanticoncentrationgoal}, under the assumption that $X_i$ have a sub-exponential distribution.  We are grateful to him for allowing us to include his elegant proof in this version of the paper.

\section{Our results}
Let us write $X=(X_{1},\ldots ,X_{n})$ so that $\<\alp,X\>=\sum_{i=1}^{n}\alp_{i}X_{i}$ and $\<\bet,X\>=\sum_{i=1}^{n}\bet_{i}X_{i}$. First we show in Section~\ref{sec:gaussian} that if $X_{i}$ are i.i.d. standard Gaussian random random variables, then
\begin{align}\label{eq:gaussianineq}
\P\{|\<\alp,X\>|\le |\<\bet,X\>|\}\le 2\frac{\|\bet\|}{\|\alp\|}.
\end{align}
This may be taken as a motivation for \eqref{eq:relanticoncentrationgoal}, but without the second term.  
As will be clear later, for discrete random variables, the second term become necessary.  Our main results are as follows:
\beit
\item When $X_i$s have a  log-concave distribution,  \eqref{eq:gaussianineq} holds with a larger  constant (Theorem~\ref{thm:logconcave}).
\item If $X_i$s have sub-Gaussian (Theorem~\ref{thm:subgaussian}) or mean zero sub-exponential distibution (Theorem~\ref{thm:subexponential}), we prove \eqref{eq:relanticoncentrationgoal}, but losing a factor of $\log(\|\alp\|/\|\bet\|)$ in the first term on the right. 
\item There are a few other minor results (Corollary~\ref{cor:mixlogconcave} and Corollary~\ref{cor:mixuniform}) got by taking mixtures of log-concave distributions etc.
\item After the first version of our paper appeared, Sasha Sodin sent us a sketch of a proof of \eqref{eq:relanticoncentrationgoal} for sub-exponential random variables. His result (Theorem~\ref{thm:sodinsargument}) improves on Theorem~\ref{thm:subexponential} by getting rid of the spurious $\log(\|\alp\|/\|\bet\|)$ factor. In some sense, this is the strongest result in this paper (except for the symmetry assumption which we were not able to get rid of).
\eeit
The Fourier analytic method of proof of Sodin is also entirely different from our other proofs. Hence we  retain Theorem~\ref{thm:subexponential} (and its short proof) and also give full details of Sodin's proof. For a reader with limited time, we recommend reading just the proofs of Theorem~\ref{thm:logconcave} and Theorem~\ref{thm:sodinsargument}.

Before stating the results, we recall the definition of  $\mb{LCD}$ as introduced by Rudelson and Vershynin. Among the minor variants of this quantity in their papers, we take the one in \cite{rudelsonvershyninicm}.

For a vector $\alp\in \R^{n}$ and a positive number $\gam$, define its {\em essential least common denominator} as 
\ba
\mb{LCD}_{\gam}(\alp)=\inf\l\{\theta>0\suchthat \mb{dist}(\theta\alp,\Z^{n})\le \min\l\{\gamma,\frac{1}{10}\|\theta\alp\|\r\}\r\}.
\ea
With this definition, Rudelson and Vershynin proved that (see Theorem~4.2 in \cite{rudelsonvershyninicm}) that if $X_{i}$ are i.i.d. random variables with $Q_{X_{1}}(1)=p<1$ and $\|\alp\|=1$, then for $S=\<\alp,X\>$, we have
\begin{align}\label{eq:rudelsonvershyninbound}
Q_{S}(\eps)\le C_{p}\l\{ \eps+\frac{1}{\mb{LCD}_{\gam}(\alp)}+e^{-c_{p}\gam^{2}}\r\}.
\end{align}
Here and elsewhere, one may make the choice $\gam\asymp \sqrt{n}$ so that the term $e^{-c\gam^{2}}$ become irrelevant (with $n$ discrete random variables, any non-trivial event will occur with at least $e^{-cn}$ probability).
\begin{theorem}\label{thm:subgaussian} Let $X_{i}$ be i.i.d. with a  sub-Gaussian distribution, i.e., $\P\{|X_{1}|\ge t\}\le Ce^{-ct^{2}}$. Assume $\E[X_{i}]=0$. Then, for any $\alp,\bet\in \R^{n}$, and any $\gam>0$, we have 
\begin{align*}
\P\{|\<\alp,X\>|\le |\<\beta,X\>|\}\le C'\l\{\frac{\|\beta\|}{\|\alpha\|}\sqrt{\log\frac{\|\alpha\|}{\|\beta\|}}+\frac{1}{\mb{LCD}_{\gam}(\alp/\|\alp\|)} + e^{-c'\gam^{2}}\r\}.
\end{align*}
where $C',c'$ depend on $C,c$.
\end{theorem}

A similar inequality holds under slightly milder conditions. A zero mean random variable $X$ is said to have {\em sub-exponential distribution} with parameters $(\nu,b)$ with $\nu>0,b>0$, if
\ba
\E[e^{\lam X}]\le e^{\lam^{2}\nu^{2}/2} \;\;\; \mb{ for }|\lam|\le \frac{1}{b}.
\ea
This is equivalent to the finiteness of the moment generating function $M(t)=\E[e^{tX}]$ for $|t|\le c$ for some $c>0$ which in turn is equivalent to exponential decay of tail probabilities $\P\{|X_1|>t\}$ (see ~\cite{bartlett} for details).
\begin{theorem}\label{thm:subexponential} Let $X_i$ be i.i.d. zero mean random variables with a sub-exponential distribution with parameters $(\nu,b)$.  Then, for any $\alpha,\beta \in \R^{n}$, and any $\gamma>0$, we have
\ba
\P\{|\<\alp,X\>|\le |\<\beta,X\>|\}\le C'\l\{\frac{\|\beta\|}{\|\alpha\|}\log\frac{\|\alpha\|}{\|\beta\|}+\frac{1}{\mbox{LCD}_{\gamma}(\alpha/\|\alpha\|)} + e^{-c'\gam^{2}} \r\}
\ea
where $C',c'$ depend on $\nu,b$.
\end{theorem}

The inequalities in these two theorems are sub-optimal, due to the presence of the logarithmic terms on the right. This comes from the fact that our proof works by separately bounding the probability that $|\<\alp,X\>|$ is small and the probability that $|\<\bet,X\>|$ is large. In case of Gaussian, or more generally log-concave densities, we are able to handle the joint distribution of $\<\alp,X\>$ and $\<\bet,X\>$ and hence the inequalities in \eqref{eq:gaussianineq} and in Theorem~\ref{thm:logconcave} below are optimal.

\begin{theorem}\label{thm:logconcave} If $X_i$ are i.i.d. with a non-degenerate  log-concave density that is symmetric about $0$, then 
\ba
\P\{|\<\alp,X\>|\le |\<\bet,X\>|\}\le C\frac{\|\bet\|}{\|\alp\|}
\ea
where $C$ is a constant.
\end{theorem}
These three theorems and Theorem~\ref{thm:sodinsargument} below are the main results of this paper. Since log-concave densities decay exponentially, in all these theorems we have exponential decay of the tails of $X_{1}$.  By taking mixtures of log-concave random variables, one can allow somewhat heavier tails, as in the following two corollaries to Theorem~\ref{thm:logconcave}.
\begin{corollary}\label{cor:mixlogconcave}
Let $X_{i}=\xi_{i} Y_{i}$ where $Y_{i}$ are i.i.d. with a symmetric, log-concave density, $\xi_{i}$ are i.i.d. positive random variables with $\E[\xi_{1}^{2}]\le B$ and $\E[1/\xi_{1}^{2}]\le B$ for some $B$ and $\xi_{i}$ are independent of $Y_{i}$s.  Then, 
\ba
\P\{|\<\alp,X\>|\le |\<\beta,X\>|\}\le CB\,\frac{\|\beta\|}{\|\alpha\|}
\ea
where $C$ is the constant in Theorem~\ref{thm:logconcave}.
\end{corollary}
In particular, writing a unimodal density as a mixture of uniform densities on intervals, we get the following conclusion. 
\begin{corollary}\label{cor:mixuniform} Let $X_{i}$ be i.i.d. with a symmetric unimodal density $f$ such that $\int t^{2}f(t)dt\le B$ and $\int_{0}^{1}t^{-3}[\P\{|X|\le t\}-2tf(t)]dt\le B$ for some $B$. Then, 
\ba
\P\{|\<\alp,X\>|\le |\<\beta,X\>|\}\le 12C B\,\frac{\|\beta\|}{\|\alpha\|}.
\ea
\end{corollary}
Note that second  condition on the density is satisfied by $f(t)=e^{-|t|^{1+\del}}$ for any $\del>0$ but not by $e^{-|t|}$. The condition restricts how sharply the density can peak at the origin.

\begin{remark} One can get a variant of Corollary~\ref{cor:mixlogconcave} with the bound of $\frac{1}{p_{\ell}}\E[(\xi_{1}^{2}+\ldots +\xi_{\ell}^{2})^{-1}]$ where $p_{\ell}$ is the $\ell$-th largest of the numbers $\alp_{i}^{2}/\sum \alp_{i}^{2}$. This is some times applicable when we have some information on $\alp$ (eg., that it is not dominated by a single $\alp_{i}$). We skip details.
\end{remark}

Now we state the result of Sodin referred to earlier. This is an improvement over Theorem~\ref{thm:subexponential}, except for the assumption of symmetry. 
\begin{theorem}\label{thm:sodinsargument} Let $X_i$ be i.i.d. zero mean random variables with a sub-exponential distribution with parameters $(\nu,b)$. Assume that the distribution of $X_i$s is symmetric about zero. Then, for any $\alpha,\beta \in \R^{n}$, and any $\gamma>0$, we have
\begin{align*}
\P\{|\<\alp,X\>|\le |\<\beta,X\>|\}\le C'\l\{\frac{\|\beta\|}{\|\alpha\|}+\frac{1}{\mbox{LCD}_{\gamma}(\alpha/\|\alpha\|)} + e^{-c'\gam^{2}} \r\}
\end{align*}
where $C',c'$ depend on the distribution of $X_1$.
\end{theorem}

\section{Proof of the inequality for Gaussians}\label{sec:gaussian}
We prove \eqref{eq:gaussianineq} in this section. Let $U'=\<\alp,X\>$ and $V'=\<\bet,X\>$. Let $\rho=\frac{\<\alp,\bet\>}{\|\alp\|\|\bet\|}$ and let $\xi,\eta$ be i.i.d. standard Gaussians. For simplicity of notation, let $\theta=\|\bet\|/\|\alp\|$.   Then $(U',V')$ has the same joint distribution as $(U,V)$ where $V=\|\bet\|\xi$ and $U=\|\alp\|(\rho \xi+\sqrt{1-\rho^{2}}\eta)$. Hence,
\begin{align*}
\P\{|U'|\le |V'|\} = \P\l\{\Big|\frac{\eta}{\xi}+\frac{\rho}{\sqrt{1-\rho^{2}}}\Big|\le \frac{\theta}{\sqrt{1-\rho^{2}}}\r\}=\P\l\{\frac{\eta}{\xi}\in [a-\ell,a+\ell]\r\}
\end{align*}
where $a=\rho/\sqrt{1-\rho^{2}}$ and $\ell=\theta/\sqrt{1-\rho^{2}}$.  Now, $\eta/\xi$ has Cauchy distribution whose density $1/\pi(1+t^{2})$ is unimodal and has the maximum value of $1/\pi$. Hence,
\begin{align}\label{eq:bdforcauchyinterval}
\P\l\{\frac{\eta}{\xi}\in [a-\ell,a+\ell]\r\}\le \begin{cases}
 \frac{2\ell}{\pi} & \mb{ for any }a,\ell, \\
 \frac{2\ell}{\pi (a-\ell)^{2}} &\mb{ if }a-\ell>0, \\
 \frac{\ell}{\pi (a+\ell)^{2}} &\mb{ if }a+\ell<0.
\end{cases}
\end{align}
If $\rho^{2}\le 1-\frac{1}{\pi^{2}}$, we use the first bound in \eqref{eq:bdforcauchyinterval}  to get
\ba
\P\{|U'|\le |V'|\} \le \frac{2\theta}{\pi\sqrt{1-\rho^{2}}}\le 2\theta.
\ea
 If $\rho^{2}>1-\frac{1}{\pi^{2}}$, then  use the second or third bound in \eqref{eq:bdforcauchyinterval} (depending on $\rho>0$ or $\rho<0$) to get
\ba
\P\{|U'|\le |V'|\}\le \frac{2\sqrt{1-\rho^{2}}}{\pi(\rho-\theta)^{2}}\theta.
\ea
We may assume $\theta\le \frac12$ (otherwise  $2\theta$ is a trivial bound for any probability). Then, checking numerically that $(\rho-\theta)^{2} \ge 0.15$ and $\sqrt{1-\rho^{2}}\le \frac{1}{\pi}$, we see that the right hand side of the previous inequality is smaller than $2\theta$.

\section{Proofs of Theorems~\ref{thm:subgaussian} and \ref{thm:subexponential}}
\bprf[Proof of Theorem~\ref{thm:subgaussian}] 
As $X_i$ are i.i.d. sub-Gaussian, by a version of Bernstein's inequality (see Theorem~3.3 in \cite{rudelsonnotes}), for any $t>0$, we have
\begin{equation}
\P\{|\<\beta,X\>|\ge t\} \le Ce^{-c\frac{t^{2}}{\|\beta\|^{2}}}.
\end{equation}
Next, using the Rudelson-Vershynin inverse Littlewood-Offord result~\eqref{eq:rudelsonvershyninbound}, we have
$$
\P\l\{|\<\alp,X\>|\leq t\r\}\leq 
c_1\l\{\frac{1}{\mbox{LCD}_{\gamma}(\alpha/\|\alpha\|)}+\frac{t}{\|\alpha\|}\r\}+c_2e^{-c_3\gamma^2}
$$
Hence,
$$
\P\{|\<\alp,X\>|\le |\<\beta,X\>|\}\leq Ce^{-c\frac{t^2}{\|\beta\|^2}}+c_1\l\{\frac{1}{\mbox{LCD}_{\gamma}(\alpha/\|\alpha\|)}+\frac{t}{\|\alpha\|}\r\}+c_2e^{-c_3\gamma^2}
$$
Choose $t=\frac{1}{\sqrt{c}}\|\beta\|\sqrt{\log\frac{\|\alpha\|}{\|\beta\|}}$, and we get 
$$
\P\{|\<\alp,X\>|\le |\<\beta,X\>|\}\lesssim \frac{\|\beta\|}{\|\alpha\|}+\frac{\|\beta\|}{\|\alpha\|}\sqrt{\log\frac{\|\alpha\|}{\|\beta\|}}+\frac{1}{\mbox{LCD}_{\gamma}(\alpha/\|\alpha\|)}+ e^{-c\gam^{2}}.
$$
We shall always take $\|\beta\|<\|\alpha\|$ so that the bound in the statement of the lemma follows.
\eprf
In proving Theorem~\ref{thm:subexponential}, we shall use the key concentration property
\ba
\P\{|X|\ge t\}\le\begin{cases}
e^{-t^{2}/2\nu^{2}} &\mb{ if }0\le t\le \nu^{2}/b,\\
e^{-t/2b} &\mb{ if }t>\nu^{2}/b.
\end{cases}
\ea
This well-known inequality (essentially due to Bernstein) may be worked out from the exercises on page 205 of Uspensky's book~\cite{uspensky}. For a  more easily accessible reference, see  \cite{bartlett}. 
For the following proof, we introduce the notation $\bet_{\max}=\max_{i\le n}|\bet_{i}|$.

\bprf[Proof of Theorem~\ref{thm:subexponential}] As $X_i$ are i.i.d. sub-exponential with parameters $(\nu,b)$, hence $\beta_iX_i$ are independent sub-exponential with parameters $(\beta_i\nu,\beta_ib)$, and $\< \beta,X\>$ is sub-exponential with parameters $(\nu^*,b^*)$ where $b^*=b\bet_{\max}$, and $\nu^*=\nu \|\beta\|$. Hence,
\begin{align}\label{eq:boundsforsumbetax}
\P\{|\<\beta,X\>|\geq t_0\}\leq 
\left\{
\begin{array}{ll}
e^{-\frac{t_0^2}{2\nu^2\|\beta\|^2}} & \mbox{ if } 0\leq t_0\leq \frac{\nu^2\|\beta\|^2}{b\bet_{max}} ,\\
e^{-\frac{t_0}{2b\bet_{\max}}} & \mbox{ if } t_0 >\frac{\nu^2\|\beta\|^2}{b\bet_{max}}.
\end{array}
\right.
\end{align}
Again, using Rudelson-Vershynin's inverse Littlewood-Offord result, we have
\begin{align}\label{eq:boundforalpXterm} 
\P\l\{\frac{|\<\alp,X\>}{\|\alpha\|}\leq \frac{t_0}{\|\alpha\|}\r\}\leq 
C_1\l\{\frac{1}{\mbox{LCD}_{\gamma}(\alpha/\|\alpha\|)}+\frac{t_0}{\|\alpha\|}\r\}+C_2e^{-c_3\gamma^2}
\end{align}
When $\frac{\bet_{\max}}{\|\beta\|}\sqrt{\log\frac{\|\alpha\|}{\|\beta\|}}\leq \frac{\nu}{\sqrt{2}b}$, put $t_0=\sqrt{2}\nu \|\beta\|\sqrt{\log\frac{\|\alpha\|}{\|\beta\|}}$ and use the first inequality in \eqref{eq:boundsforsumbetax}. That term become $\|\bet\|/\|\alp\|$. Adding it to \eqref{eq:boundforalpXterm} gives us the bound
\ba
\P\{|\<\alp,X\>|\le |\<\bet,X\>|\}&\le \frac{\|\bet\|}{\|\alp\|} + C_1\l\{\frac{1}{\mbox{LCD}_{\gamma}(\alpha/\|\alpha\|)}+\sqrt{2}\nu\frac{\|\bet\|}{\|\alpha\|}\sqrt{\log\frac{\|\alpha\|}{\|\beta\|}}\r\}+C_2e^{-c_3\gamma^2}
\\
&\le  C'\l\{\frac{\|\beta\|}{\|\alpha\|}\sqrt{\log\frac{\|\alpha\|}{\|\beta\|}}+\frac{1}{\mbox{LCD}_{\gamma}(\alpha/\|\alpha\|)} + e^{-c'\gam^{2}} \r\}
\ea
which is better than we claimed, because of the square root on the logarithmic factor.   When $\frac{\bet_{\max}}{\|\beta\|}\sqrt{\log\frac{\|\alpha\|}{\|\beta\|}}> \frac{\nu}{\sqrt{2}b}$, put $t_0=2b\bet_{\max} \log \frac{\|\alpha\|}{\|\beta\|}$ and use the second inequality in \eqref{eq:boundsforsumbetax}. That term is again $\|\bet\|/\|\alp\|$. Adding it to \eqref{eq:boundforalpXterm} gives us the bound
\ba
\P\{|\<\alp,X\>| \le |\<\bet,X\>|\} &\le \frac{\|\bet\|}{\|\alp\|}+C_1\l\{\frac{1}{\mbox{LCD}_{\gamma}(\alpha/\|\alpha\|)}+2b\frac{\bet_{\max}}{\|\alpha\|}\log \frac{\|\alpha\|}{\|\beta\|}\r\}+C_2e^{-c_3\gamma^2}\\
&\le C'\l\{\frac{\|\bet\|}{\|\alpha\|}\log\frac{\|\alpha\|}{\|\beta\|}+\frac{1}{\mbox{LCD}_{\gamma}(\alpha/\|\alpha\|)} + e^{-c'\gam^{2}} \r\}.
\ea
since $\bet_{\max}\le \|\bet\|$. This completes the proof.
\eprf

\section{Proof of Theorem~\ref{thm:logconcave}}
A probability distribution $\mu$ on $\R^{2}$ is said to be isotropic if it has zero mean and identity covariance, ie.,
\ba
\int_{R^{2}}\x d\mu(\x)=0, \mb{ and }\int_{\R^{2}} \x \x^{t}d\mu(\x) = \mat{1}{0}{0}{1}.
\ea
We shall use the following lemma about isotropic log-concave measures in the plane. Let $D(x,r)$ denote the open disk of radius $r$ centered at $x$.
\begin{lemma}\label{lem:bivariatelogconcave}
 Let $p(x,y)=e^{-f(x,y)}$ be an isotropic, log-concave density on $\R^{2}$.  Let $\mathcal{L}=\{(x,y) : p(x,y)\geq p(0,0)/2\}$. 
There exist two numerical constants $0<a<A$ and $0<b<B$, such that $D(0,a)\subseteq \mathcal{L} \subseteq D(0,A)$ and $b\le \max\limits_{(u,v)}p(u,v)\le B$.
\end{lemma}
\bprf This can be read off from Lemma~5.14 of Lovasz and Vempala~\cite{lovaszvempala} (their lemma is valid in any dimension) as follows: Part (a) of that lemma immediately gives $a=1/9$. Next, by part (d) of their lemma, $p(0,0)\ge 2^{-14}$. Integrating the density over $\mathcal L$, we see that $\mb{area}(\mathcal L)\le 2^{15}$. If $\mathcal L$ intersects $\partial D(0,A)$ at a point $\x$, then by convexity (draw the tangents from $\x$ to the circle $\partial D(0,a)$ and join $\x$ and these points of tangency to  the origin to get two right angles triangles) its area is at least $a\sqrt{A^{2}-a^{2}}$. Hence, we must have $A\le a^{-1}2^{16}$. 

Lastly, by the already quoted bound, we may take $b=2^{-14}$ and $B=2/\pi a^{2}$ (the latter because  the density is at least $p(0,0)/2$ on $D(0,a)$). 
\eprf
\para{Sketch of an alternate argument} If one does not care about explicit constants, it is also possible to prove Lemma~\ref{lem:bivariatelogconcave} by a compactness argument. We explain it to show the existence of the number $a>0$. It is clear that for any isotropic, log-concave density, there is an $a>0$ that works, what is non-trivial is the uniform choice of the constant. Now suppose there is no such uniform constant $a>0$. Then we may take a sequence of isotropic, log-concave densities $p_{n}$ such that $p_{n}(x_{n})\le  p_{n}(0,0)/2$ with $|x_{n}|\le 1/n$. Since rotation of an isotropic log-concave density is also isotropic and log-concave, we may assume that $x_{n}=(1/n,0)$. The space of log-concave measures is closed under weak convergence (Proposition~3.6 of \cite{saumardwellner}), hence we may assume that $p_{n}(x)dx$ converge weakly to a log-concave measure $\mu$. For log-concave measures, weak convergence implies convergence of all moments (Corollary~6 in the arXiv version of \cite{meckeses}), hence   $\mu$ is isotropic. But now, the density of $\mu$ must vanish on the $(0,\infty)$, which contradicts the existence of $a$ specific to $\mu$.  This shows the existence of a uniform constant $a>0$ as claimed. Similarly one can argue for the existence of $A$, $b$ and $B$.

Now we turn to the proof of Theorem~\ref{thm:logconcave}.

\para{Claim} If Theorem~\ref{thm:logconcave} holds when $\<\alp,\bet\>=0$, then it hold for any $\alp,\bet$.

\bprf Given any $\alp,\bet$ (not necessarily orthogonal), write  $\beta=a\alpha + \gamma$ where $\<\alpha,\gamma\>=0$. Since $|\langle \beta,x \rangle| \leq |a||\langle \alpha,x\rangle| + |\langle \gamma,x \rangle |$, we get
 $$\P(|\langle\alpha,x\rangle| \leq |\langle\beta,x\rangle|) \leq 
 \P \left(|\langle\alpha,x\rangle|\leq \frac{|\langle\gamma,x\rangle|}{1-|a|}\right) \leq C\frac{ \|\gamma\|}{\|\alpha\|(1-|a|)}$$ 
 where the last inequality holds because  $\langle\alpha,\frac{\gamma}{1-|a|}\rangle=0$, and our assumption that relative LO holds when inner product is $0$.

And now as $\|\beta\|^2=|a|^2\|\alpha\|^2 +\|\gamma\|^2$, hence $\|\gamma\| \leq \|\beta\|$, and
$|a| \leq \|\beta\|/\|\alpha\| <1/10$ ( without loss of generality we can assume this, otherwise we can take the constant C in the RHS of the relative LO inequality to be greater than 10, so that the RHS becomes greater than 1, and hence the inequality holds trivially). Hence 
\ba
C \frac{\|\gamma\|}{\|\alpha\|(1-|a|)} \leq C'\frac{ \|\beta\|}{\|\alpha\|}.
\ea
Thus, it suffices to prove  Theorem~\ref{thm:logconcave} when $\langle\alpha,\beta\rangle=0$.
\eprf
Now we prove the theorem for orthogonal $\alp,\bet$.
\begin{proof}[Proof of Theorem~\ref{thm:logconcave} when $\<\alp,\bet\>=0$]
If $\alp,\bet$ are orthogonal and non-zero vectors, then define $U=\<\alp,X\>/\|\alp\|$ and $V=\<\bet,X\>/\|\bet\|$. Clearly $(U,V)$ has an isotropic, log-concave distribution. Hence,
\ba
\P\{|\<\alp,X\>|\le |\<\bet,X\>|\}=\P\{(U,V)\in S\}
\ea
where $S=\{(u,v)\suchthat \frac{|u|}{|v|}\le \|\bet\|/\|\alp\|\}$. Note that $S$ is a union of two sectors in the plane, each with an angle of $2\theta$ where $\tan \theta=\|\bet\|/\|\alp\|$.  By Lemma~\ref{lem:bivariatelogconcave} and the log-concavity of $p$, we have the bound $p(u,v)\le p(0,0)2^{-k}\le B2^{-k}$ for $(u,v)\not\in D(0,kA)$ and $k\ge 1$. On $D(0,A)$ we use the bound $p(u,v)\le B$. Hence, 
\ba
\P\{(U,V)\in S\}&\le \sum_{k=1}^{\infty} B2^{-k}\mb{area}(S\cap D(0,kA))\\
&= 2\pi \theta BA^{2}\sum_{k=1}^{\infty}k^{2}2^{-k} \\
&\le C\theta
\ea
for some $C$. As $\theta \le \tan\theta$ and $\tan\theta=\frac{\|\bet\|}{\|\alp\|}$, we get $\P\{(U,V)\in S\}\le C\frac{\|\bet\|}{\|\alp\|}$.
\eprf

\section{Proofs of Corollary~\ref{cor:mixlogconcave} and Corollary~\ref{cor:mixuniform}}
\bprf[Proof of Corollary~\ref{cor:mixlogconcave}] Write $X_{i}=\xi_{i}Y_{i}$ where $Y_{i}$ are i.i.d. with a log-concave distribution. Condition on $\xi_{i}$s and apply Theorem~\ref{thm:logconcave} to get
\begin{align}\label{eq:conditionandtakeexpt}
\P\l\{|\<\alp,X\>|\le |\<\bet,X\>|\r\} &= 10\E\l[\sqrt{\frac{\sum_{i=1}^{n}\bet_{i}^{2}\xi_{i}^{2}}{\sum_{i=1}^{n} \alp_{i}^{2}\xi_{i}^{2}}} \r]\le 10\sqrt{\E\l[\sum_{i=1}^{n}\bet_{i}^{2}\xi_{i}^{2}\r]}\sqrt{\E\l[\frac{1}{\sum_{i=1}^{n}\alp_{i}^{2}\xi_{i}^{2}}\r]}
\end{align}
by Cauchy-Schwarz inequality. Now, by the bound $\E[\xi_{i}^{2}]\le B$, we get
\ba
\E\l[\sum_{i=1}^{n}\bet_{i}^{2}\xi_{i}^{2}\r] \;\le\; B\sum_{i=1}^{n}\bet_{i}^{2} \; = \; B\|\bet\|^{2}.
\ea
By Jensen's inequality aplied to the convex function $x\mapsto 1/x$, we get
\ba
\E\l[\frac{1}{\sum_{i=1}^{n}\alp_{i}^{2}\xi_{i}^{2}}\r] \; \le  \; \frac{1}{\|\alp\|^{2}}\sum_{i=1}^{n}\frac{\alp_{i}^{2}}{\|\alp\|^{2}} \E[1/\xi_{i}^{2}] \; \le  \; \frac{B}{\|\alp\|^{2}}.
\ea
Using these bounds, we see that the right hand side of \eqref{eq:conditionandtakeexpt} is at most $10B\|\bet\|/\|\alp\|$.
\eprf

\bprf[Proof of Corollary~\ref{cor:mixuniform}] For $0<y<f(0)$, let $g(y)$ be the length of the interval $\{s\suchthat f(s)\ge y\}$. Then, $g$ is a density (evaluate the area under $f$ by integrating over the x-coordinate first and then over the y-coordinate). Further, if $\xi$ is a random variable with density $g$ and $Y$ has $\mb{Uniform}[-1/2,1/2]$ density and $\xi$, $Y$ are independent, then $\xi Y$ has the density $f$. Since $Y$ is log-concave, we can apply Corollary~\ref{cor:mixlogconcave} to get the conclusion we want, 
 {\em if} $\E[\xi^{2}]$ and $\E[\xi^{-2}]$ are finite.

Since $\E[Y^{2}]=1/12$, we see that $\E[\xi^{2}]=12\int t^{2}f(t)dt$. Further, 
\ba
\E[\xi^{-2}]=\int_{0}^{\infty}\P\{\xi<t\}\frac{2}{t^{3}}dt = \int_{0}^{\infty}\frac{2}{t^{3}}[\P\{|X|\le t\}-tf(t)] dt.
\ea
Hence the conditions in the statement of the theorem ensure that $\E[\xi^{2}]$ and $\E[\xi^{-2}]$ are finite, and the conclusion follows.
\eprf

\section{Sodin's proof of Theorem~\ref{thm:sodinsargument}}
By scaling $\alp$ and $\bet$ to have unit norm, we recast the  theorem in the following equivalent form: Let $U=\<\alp,X\>$ and $V=\<\bet,X\>$ where $\|\alp\|=\|\bet\|=1$. Then there are constants $C,c$ depending on the distribution of $X_1$ such that for any $\gam>0$, we have
\begin{align}\label{eq:equivalentformoftheorem7}
\P\l\{|U|\le \eps |V|\r\}\le  C\l\{\eps +\frac{1}{\mb{LCD}_{\gam}(\alp)}+e^{-c\gam^2}\r\}\;\;\; \mb{ for any }\eps>0.
\end{align}
We may also replace $X_i$ by $X_i/b$ and assume that they are sub-exponential with parameters $(\nu,1)$. Thus if $\phi(\lam)=\E[e^{i\lam X_1}]$ denotes the characteristic function of $X_i$s and $M(\lam)=\phi(-i\lam)$ denotes the moment generating function, then  $M(\lam)\le e^{\lam^2\nu^2/2}$ for $|\lam|\le 1$. As stated in \eqref{eq:boundsforsumbetax}, this implies that 
\begin{align}\label{eq:tailboundforsumofsubexprvs}
\P\{|V|>u\}\le \begin{cases} e^{-\frac{u^2}{2\nu^2}} & \mb{ for } 0\le u\le \nu^2, \\
e^{-\frac{u}{2}} & \mb{ for  }u>\nu^2.\end{cases}
\end{align}
since $\beta_{\max}\le \|\bet\|=1$.

Fix $\eps>0$ and break the event in \eqref{eq:equivalentformoftheorem7} as follows.
\begin{align*}
\P\{|U|<\eps |V|\}\le \P\{|U|<\eps\}+\sum_{k=0}^{\infty} \P\l\{|U|<2^{k+1}\eps, \ 2^k\le |V|\le 2^{k+1}\r\}.
\end{align*}
%
By the Rudelson-Vershynin inquality \eqref{eq:rudelsonvershyninbound}, the  first event can be controlled as
\begin{align}\label{eq:boundforsmallU}
\P\{|U|<\eps\}\le C\l\{ \eps+\frac{1}{L_{\gam}(\alp)}+e^{-c\gam^{2}}\r\}.
\end{align}
where we have written $L_{\gam}$ for $\mb{LCD}_{\gam}(\alp)$, for simplicity of notation. 
 We claim that for any $R\ge 1$
\begin{align}\label{eq:boundforintermediateUV}
\P\l\{|U|\le \eps R, \ V>R\r\}\le C  e^{-R}\l\{\eps R + \frac{1}{L_{\gam}} + e^{-c\gam^2}\r\}.
\end{align}
Identical bound holds for $\P\{|U|\le \eps R, \  V<-R\}$ by symmetry. Summing these estimates over $R=2^k$ (and changing $\eps$ to $2\eps$) we get
\begin{align*}
\sum_{k=0}^{\infty} \P\l\{|U|<2^{k+1}\eps, \ 2^k\le |V|\le 2^{k+1}\r\} \le  C\l\{\eps+\frac{1}{L_{\gam}}+e^{-c\gam^2}\r\}.
\end{align*}
Adding this to \eqref{eq:boundforsmallU}, we get \eqref{eq:equivalentformoftheorem7}. Thus, only the proof of \eqref{eq:boundforintermediateUV} remains.

\noindent{\em Proof of \eqref{eq:boundforintermediateUV}:} If $|U|<\eps R$ and $V>R$, then $V-\frac{1}{(\eps R)^2}U^2\ge R-1$. Therefore,
\begin{align*}
\P\l\{|U|<\eps R , V>R\r\}&\le e^{1-R}\E\l[ e^{V-\frac{1}{(\eps R)^2}U^2}\r] \\
&=\frac{1}{ \sqrt{\pi}}e^{1-R} \int_{\R}\E\l[ e^{V+2i\frac{x}{\eps R}U}\r] e^{-x^2}dx
\end{align*}
using the identity $\int_{\R}e^{2itx- x^2}dx=\sqrt{\pi}e^{-t^2}$ and interchanging the integral and expectation. Write $V+2i\frac{x}{\eps R}U = \sum_{k=1}^n \l(\bet_k+2i\frac{x\alp_k}{\eps R}\r)X_k$ to see that
\begin{align*}
\E\l[ e^{V+2i\frac{x}{\eps R}U}\r] &=\prodd_{k=1}^n\phi\l(-i\bet_k+2\frac{\alp_k x}{\eps R}\r) \\
&= \prodd_{k=1}^n M(\bet_k) \; \prodd_{k=1}^n\phi_k\l(\frac{2\alp_k x}{\eps R}\r)
\end{align*}
where $\phi_k(\lam):=\frac{1}{M(\bet_k)}\phi\l(\lam-i\bet_k\r)$ is the characteristic function of the exponentially tilted measure $dF_k(x):=\frac{1}{M(\bet_k)}e^{\bet_k x }dF(x)$, with $F$ being the distribution of $X_1$. As $|\bet_k| \le 1$,  we have $M(\beta_k)\le e^{\nu^2\bet_k^2/2}$ for each $k$. Using $\|\bet\|=1$ the product of $M(\bet_k)$ over $k$ is at most $e^{\nu^2/2}$.
Consequently, writing $t_k=2\alp_k/\eps R$,
\begin{align}\label{eq:boundintermsofintegralofproductofcfs}
\P\l\{|U|<\eps R , V>R\r\}&\le C e^{-R}\intt_{\R} \prodd_{k=1}^n \l|\phi_k\l(t_k x\r)\r| \; e^{-x^2}dx.
\end{align}
We introduce some notation. Let $Y_k,Y_k'$ denote independent random variable with distribution $F_k$ and let $W_k=Y_k-Y_k'$. Then $|\phi_k|^2$ is the characteristic function  of $W_k$. Fix $\del>0$ and $p<1$ such that $Q_{X_1}(\del)\le p$.  Let $q_k=\P\{|W_k|\ge 2\del\}$.  Then
\begin{align*}
\log |\phi_k(t)|^2 &\le - (1-|\phi_k(t)|^2) \\
&=-\E[1-\cos(tW_k)]  \\
&\le -q_k\E\l[1-\cos(tW_k)\Given |W_k|\ge 2\del\r].
\end{align*}
By Lemma~\ref{lem:tiltedanduntilted} and its Corollary~\ref{cor:tiltedanduntilted} that are proved later, using the bound $M(t)\le M(1)$ for $|t|\le 1$, we deduce that there are positive constants $q$ and $\tau$ depending only on $F$ such that for all $k$ and for all $s$ we have
\begin{align*}
q_k\ge q \;\;\mb{ and }\;\; \E[1-\cos(sW_k)\Given |W_k|\ge 1]\ge \tau \E[1-\cos(sW)\Given |W|\ge 1],
\end{align*}
 where $W=X_1-X_1'$ (the analogue of $W_k$ but for the untilted random variable).  Using these uniform estimates in  \eqref{eq:boundintermsofintegralofproductofcfs}, we arrive at
\begin{align*}
\P\l\{|U|<\eps R , V>R\r\}&\le C e^{-R}\int_{\R}\exp\l\{-\frac{q\tau}{2} \E\l[ \summ_{k=1}^n(1-\cos(t_kxW)) \given  |W|\ge \del\r] \r\} e^{-x^2}\ dx \\
&\le C e^{-R}\int_{\R} \E\l[\exp\l\{-\frac{q\tau}{2} \summ_{k=1}^n(1-\cos(t_kxW))\r\}   \given |W|\ge \del \r]e^{-x^2}dx
\end{align*}
by Jensen's inequality. Now interchange conditional expectation with integral and then replace the conditional expectation over $|W|\ge \del$ by the maximum over  $|W|\ge \del$. That gives us 
\begin{align}\label{eq:boundintermsofcosines}
\P\l\{|U|<\eps R , V>R\r\} &\le  C e^{-R}\sup_{|w|\ge \del}  \intt_{\R} \exp\l\{-\frac{q\tau}{2} \summ_{k=1}^n(1-\cos(t_kxw))    -x^2\r\} dx
\end{align}
From this point, the arguments are virtually identical to those of Friedland and Sodin~\cite{friedlandsodin} (one small difference is that their version of LCD is not the same). Since $1-\cos(\theta)\ge 8\ \mb{dist}^2(\frac{\theta}{2\pi},\Z)$, we have $ \sum_{k=1}^n(1-\cos(t_kxw)) \ge 8 \ \mb{dist}^2\l(\frac{xw}{\pi \eps R}\alp,\Z^n\r)$. Fix $w\ge \del$ (identical argument applies to $w\le -\del$) and use this bound in the integral above to write
\begin{align} \label{eq:boundonintegralforfixedw}
 \intt_{\R} \exp\l\{-\frac{q\tau}{2} \summ_{k=1}^n(1-\cos(t_kxw))    -x^2\r\} dx &\le  \intt_{\R} \exp\l\{-4 q\tau  \mb{dist}^2\l(\frac{xw}{\pi \eps R}\alp,\Z^n\r)\r\}  e^{ -x^2} dx \nonumber \\
  =8q\tau & \int_{0}^{\infty}\mu\l\{x\suchthat  \mb{dist}\l(\frac{xw}{\pi \eps R}\alp,\Z^n\r)\le z\r\} \ z e^{-4q\tau z^2} \; dz
\end{align}
where $\mu$ is the measure $e^{-x^2}dx$ on the line (the last equality is by  the well-known principle $\int f d\mu=\int_0^{\infty}\mu\{f>t\}dt$ for non-negative $f$).

Let $I(z):=\{x\in \R\suchthat  \mb{dist}\l(\frac{xw}{\pi \eps R}\alp,\Z^n\r)\le z\}$. 
 For $z\le \frac12 \gam$, we now show that $I(z)$ is a union of well-separated short intervals. Indeed, if $x,y\in I(z)$, then $\frac{(x-y)w}{\pi \eps R}\alp$ is within $2z$ distance of $\Z^n$. Hence, by the definition of $\mb{LCD}$, we must have 
\begin{align*}
\mb{ either } \;\;\; \frac{|x-y|w}{\pi \eps R}\ge L_{2z}\ge  L_{\gam} \;\;\; \mb{ or } \;\;\; \frac{1}{10}\frac{|x-y|w}{\pi \eps R}\le 2z.
\end{align*}
Therefore, $I(z)$ is contained in a union of intervals $I_j$, $j\in \Z$, such that \begin{inparaenum}[(a)]\item $I_j$ lies to the left of $I_{j+1}$, \item  each $I_j$ has length at most $\frac{20 \pi \eps R z}{w}$ and \item $I_j$ and $I_{j+1}$ are at distance at least $\frac{\pi \eps R}{w}L_{\gam}$ from each other. \end{inparaenum} Indexing them so that $a_0$ is the closest among $a_j$s to the origin, we see that $I_j$ is at a distance of at least $(|j|-1)\frac{\pi \eps R}{w}L_{\gam}$ from the origin. Thus,
\begin{align*}
\mu\{I(z)\} &\le \sum_{j=0}^{\infty} \intt_{I_j}e^{-x^2}dx  \\
&\le \frac{20\pi \eps R z}{w}\l(1+2\sum_{j=1}^{\infty} \exp\l\{-\frac{1}{w^2}(|j|-1)^2 \pi^2 \eps^2 R^2 L_{\gam}^2 \r\}\r)
\end{align*} 
By a standard comparison of the sum to the integral, we get
\begin{align*}
\mu\{I(z)\} &\le \frac{20\pi \eps R z}{w} \l(1+2\int_0^{\infty} \exp\l\{-\frac{1}{w^2}u^2 \pi^2 \eps^2 R^2 L_{\gam}^2\r\} du \r) \\
&=  \frac{20\pi \eps R z}{w} \l(1+\frac{\sqrt{\pi} w}{ \pi \eps R L_{\gam}} \r) \\
&\le 70 z\l(\frac{\eps R}{w} +\frac{1}{L_{\gam}} \r).
\end{align*}
Plugging this bound (and the trivial bound $\mu\{I(z)\}\le \mu(\R) = \sqrt{\pi}$ for $z\ge \frac12 \gam$) into \eqref{eq:boundonintegralforfixedw} to bound that integral as
\begin{align*}
&\le 70 \l(\frac{\eps R}{w} +\frac{1}{L_{\gam}} \r)\intt_{0}^{\gam/2}8q\tau z^2 e^{-4q\tau z^2} dz +\int_{\gam/2}^{\infty} 8q\tau z e^{-4q\tau z^2} dz & \\
&\le C \l(\frac{\eps R}{w}+\frac{1}{L_{\gam}}+e^{-c\gam^2}\r)
\end{align*}
where $C,c$ depend on $q$ and $\tau$. Since $w\ge \del$, absorbing $1/\del$ into $C$, from \eqref{eq:equivalentformoftheorem7} we have
\begin{align*}
\P\{|U|<\eps R, V>R\}\le C e^{-R}\l(\eps R +\frac{1}{L_{\gam}}+e^{-c\gam^2}\r)
\end{align*}
 This completes the proof of \eqref{eq:boundforintermediateUV}. \hfill \qed

The following lemma and its corollary were used in the proof.  Its content is that the exponential tilts of a given probability distribution are uniformly comparable to the original distribution, as long as the tilting parameter is bounded. We assume  symmetry here (by a variant of this Lemma without symmetry, one may enable one to remove the symmetry assumption in Theorem~\ref{thm:sodinsargument}, but we do not know how).
\begin{lemma}\label{lem:tiltedanduntilted} Let $F$ be a  probability distribution on the line symmetric about $0$. Let $dF_t(x)=\frac{1}{M(t)}e^{tx}dF(x)$ where $M(t)=\int e^{tx}dF(x)$.  Let $\phi:\R\mapsto [0,1]$ be an even measurable function. Then, for any $t\in \R$
\begin{align*}
\int\!\!\int \phi(x-x')\ dF(x)dF(x') \le M(t)\int\!\!\!\int \phi(x-x') \ dF_t(x)dF_t(x').
\end{align*}
\end{lemma}
\bprf Write
\begin{align*}
&\int\!\!\!\int \phi(x-x')\ dF(x)dF(x') =    \int\!\!\int \phi(x-x')e^{\frac12 t(x+x')} e^{-\frac12 t(x+x')}dF(x)dF(x') \\
& \;\;\; \le \l(\int\!\!\!\int \phi(x-x')e^{ t(x+x')}dF(x)dF(x') \r)^{\frac12} \l(\int\!\!\!\int \phi(x-x')e^{-t(x+x')}dF(x)dF(x') \r)^{\frac12}.
\end{align*}
If we make the change of variables $(x,x')\mapsto (-x,-x')$ in the second integral, then the evenness of $\phi$ and the symmetry of $F$ shows that it is identical to the first integral. Thus the right hand side is equal to 
$M(t)\int\!\!\!\int \phi(x-x') \ dF_t(x)dF_t(x')$. 
\eprf
\begin{corollary}\label{cor:tiltedanduntilted} In the setting of Lemma~\ref{lem:tiltedanduntilted}, let $X_t,X_t'$ be i.i.d. random variables with distribution $F_t$. Let $W_t=X_t-X_t'$. Then $\P\{|W_t|\ge \del\}\ge \frac{1}{M(t)}\P\{|W|\ge \del\}$ for some $c>0$ and all $\del>0$. 
\end{corollary}
\bprf Take $\phi(w)=\one_{|w|>\del}$ in the Lemma.
\eprf

\end{document}